
\documentstyle{amsppt}
\vsize=7.5in
\magnification=\magstep1
\def\bc{{\Bbb C}} 
\def\bp{{\Bbb P}}

\def\dd{{\Cal D}}

\def\hom{\operatorname{Hom}}

\def\mf{{meromorphic function }}
\def\mff{{f=\frac{P}{Q}}}
\def\hip{{\bc\bp^n}}
\def\hi{{\bc\bp^1}}

\newcount\parno \parno=1
\newcount\prono \prono=1
\def\sec{\S\the\parno.-\ \global\prono=1}
\def\etiqueta{\hbox{(\the\parno.\the\prono)}}
\def\finparrafo{\global\advance\parno by1
\vskip.1truecm\ignorespaces}
\def\finparrafo{\global\advance\parno by1
\vskip.1truecm\ignorespaces}
\def\cita{\ignorespaces\ \the\parno.\the\prono%
\global\advance\prono by 1}

\NoBlackBoxes

\document

\NoBlackBoxes
\topmatter

\title
On atypical values and local monodromies of meromorphic functions
\endtitle
\rightheadtext{Local monodromies of meromorphic functions}
\author
S.M. Gusein--Zade, I. Luengo, A. Melle--Hern\'andez
\endauthor

\address
Moscow State University,
Faculty of Mathematics and Mechanics,
Moscow, 119899, Russia.
\endaddress
\email
sabir\@ium.ips.ras.ru
\endemail
\address
Departamento de \'Algebra,
Facultad de Matem\'aticas,
Universidad Complutense.
E-28040 Madrid, Spain.
\endaddress
\email
iluengo\@eucmos.sim.ucm.es
\endemail
\address
Departamento de Geometr\'{\i}a y Topolog\'{\i}a,
Facultad de Matem\'aticas,
Universidad Complutense.
E-28040 Madrid, Spain.
\endaddress
\email
amelle\@eucmos.sim.ucm.es
\endemail

\thanks
First author was partially supported by Iberdrola, INTAS--96--0713, RFBR 96--15--96043.
Last two authors were partially supported by DGCYT PB94-0291.
\endthanks

\abstract
A meromorphic function on a compact complex analytic manifold
defines a $C^\infty$ locally trivial fibration over the
complement of a finite set in the projective line $\hi$. We describe
zeta-functions of local monodromies of this fibration around atypical
values.  Some applications to polynomial functions on $\bc^n$ are described.
\endabstract
\keywords Meromorphic function, monodromy, zeta-function
\endkeywords
\endtopmatter

\document

\head\sec Introduction
\endhead

We want to consider fibrations defined by meromorphic functions. In order to
have more general statements we prefer to use the notion of a \mf slightly
different from the standard one. Let $M$ be an $n$-dimensional compact complex
analytic manifold.

\definition{Definition} A {\it \mf} $f$ on the manifold $M$ is a ratio
$\frac{P}{Q}$ of two non-zero sections of a line bundle $\Cal L$ over $M.$
Meromorphic functions $\mff$ and $f'=\frac{P'}{Q'}$ (where $P'$ and $Q'$ are
sections of a line bundle $\Cal L'$) are {\it equal} if $P=UP'$ and $Q=UQ'$
where $U$ is a section of the bundle $\hom ({\Cal L'},{\Cal L})={\Cal L}\otimes
{\Cal L'}^*$ without zeroes (in particular, this implies that the bundles
${\Cal L}$ and ${\Cal L'}$ are isomorphic).
\enddefinition

A particular important case of meromorphic functions are rational functions
$\frac{P(x_0,\ldots,x_n)}{Q(x_0,\ldots,x_n)}$ on the projective space $\hip$
($P$ and $Q$ are homogeneous polynomials of the same degree).

A \mf $\mff$ defines a map $f$ from the complement
$M\setminus\{P=Q=0\}$ of the
set of common zeroes of $P$ and $Q$ to the complex line $\hi.$ The
indeterminacy set
$\{P=Q=0\}$ may have components of codimension one. For $c\in\hi$, let
$F_c=f^{-1}(c)$.

The standard arguments (using a resolution of singularities; see, e.g.,
\cite{7}) give the following statement.

\proclaim{Theorem 1} The map $f:M\setminus\{P=Q=0\}\to\hi$ is a
$C^\infty$ locally trivial fibration outside
 a finite subset of the projective line $\hi$.
 \endproclaim

Any fibre $F_{gen}=f^{-1}(c_{gen})$ of this fibration is called a generic fibre of the
\mf $f.$
The smallest subset $B(f)\subset\hi$ for which $f$ is a $C^\infty$ locally
trivial fibration over $\hi\setminus B(f)$ is called the {\it
bifurcation set} of the \mf $f.$ Its elements are called {\it atypical values.}

\remark{Remark 1} In addition to some other advantages (see Statement 1
and Remark 4) the described definition of a \mf permits to treat in the
same way the following situation. Let $\mff$ be a \mf on $M$ and let
$H=\{R=0\}$ be a hypersurface in
$M$ ($R$ is a section of a line bundle).
One can be interested to study the map defined by the restriction of the \mf $f$
to
$M\setminus (H\cup\{P=Q=0\}).$ Substituting the
\mf $f$ by $f'=\frac{P\,R}{Q\,R}$ one reduces the situation to the discussed
one, i.e., the indeterminacy set of the meromorphic function $f'$
coincides with
$H\cup\{P=Q=0\}$ and the meromorphic functions $f$ and $f'$ coincide outside it.
\endremark
\vskip 6pt

The fundamental group $\pi_1(\hi\setminus B(f))$ of the complement to the
bifurcation set acts on homology groups $H_*(F_{gen};\bc)$ of the generic fibre of the
meromorphic function $f$. The image of the group $\pi_1(\hi\setminus B(f))$ in the group of
automorphisms of  $H_*(F_{gen};\bc)$ is called the {\it monodromy group} of the
meromorphic function $f.$ It is generated by local monodromies corresponding to
atypical values (see \cite{2}).

\finparrafo

\head\sec Zeta-functions of local monodromies
\endhead

For a map $h:X\to X$ of a topological space $X$ (say, with finite dimensional
homologies) into itself, its zeta-function
$\zeta_h(t)$ is the rational function defined by
$$\zeta_h(t)=
\prod_{q\geq 0} \{ \det\,[\,id-t\,h_{*}|_{H_q(X;\bc)}]\}^{(-1)^q}.$$

\remark{Remark 2} This is the definition of the zeta-function from
\cite{2}. The zeta-function defined in \cite{1} is the inverse of this one.
\endremark
\vskip 6pt

 Let
$\zeta_f^c(t)$ be the zeta-function of the local monodromy corresponding
to the value
$c\in\hi$ (i.e., defined by a simple loop around
$c$).

\remark{Remark 3} Local monodromy and the corresponding
zeta-function are defined for any value $c\in\hi,$ not only for atypical
ones.
For a typical value of the \mf $f$,  the local monodromy is the identity
and its zeta-function is equal to $(1-t)^{\chi(F_{gen})}.$
\endremark
\vskip 6pt

The following statement is a direct consequence of the definitions.

\proclaim{Statement 1} {\rm Let $\pi:{\widetilde M}\to M$ be an analytic map of an
$n$-dimensional compact complex manifold $\widetilde M$ which is an isomorphism
outside of the union of the indeterminacy set  $\{P=Q=0\}$
of the \mf
$f$ and of a finite number of level sets $f^{-1}(c_i).$ Let ${\widetilde
f}=\frac{P\circ\pi}{Q\circ\pi}$ be the lifting of the \mf $\mff$ to
$\widetilde M$. Then the generic fibre of ${\widetilde f}$ coincides with
that of
$f$ and for each
$c\in\hi$ one has
$$\zeta_{\widetilde f}^c(t)=\zeta_f^c(t).$$
\endproclaim

\remark{Remark 4} Even if the indeterminacy set
$\{P=Q=0\}$ of the \mf
$f$ has codimension two (i.e., if the hypersurfaces $\{P=0\}$ and $\{Q=0\}$
have no common components), in general, this is not the case for the lifting
$\widetilde f.$ This is a reason for our definition of a meromorphic
function. If one starts from the usual definition, the lifting $\widetilde f$
of the \mf $f$ can be defined at some points of $\pi^{-1}(\{P=Q=0\}).$
In this case a generic level set of the meromorphic function $\widetilde f$ differs from that of $f$ and the
Statement 1 does not hold. The simplest example is $f=\frac{x}{y}$ in affine
coordinates on $\bc\bp^2,$ $\pi$ is the blowing-up of the origin in this affine
chart.
\endremark
\vskip 6pt

In order to have somewhat more attractive and unified formulae we would like to
use the notion of the integral with respect to the Euler characteristic
(\cite{8}).Let $A$ be an abelian group with the group operation $*$,
let $X$ be
a semianalytic subset of a complex manifold. Let $\Psi:X\to A$ be a
function on $X$ with values in
$A$ for which there exists a finite partitioning of $X$ into semianalytic sets
(strata) $\Xi$ such that $\Psi$ is constant on each stratum $\Xi$ (and equal to
$\psi_\Xi$). Then by definition
$$\int_X \Psi(x)\,d\chi=\sum_{\Xi}\chi(\Xi)\psi_\Xi,$$
where $\chi(\Xi)$ is the Euler characteristic of the stratum $\Xi.$
In the above formula we use the additive notations for the
operation $*.$
In what follows this definition will be used for integer valued functions and
also for local zeta-functions $\zeta_x(t)$ which are elements of the abelian
group of non-zero rational functions in the variable $t$ with respect to
multiplication. In this case in the multiplicative notations the above
formula means
$$\int_X \zeta_x(t)\,d\chi=\prod_{\Xi}(\zeta_\Xi(t))^{\chi(\Xi)}.$$

In \cite{5}, for a germ of a \mf $\varphi=\frac{F}{G}$ on $(\bc^n,0)$, there were
defined two Milnor fibres (the zero and the infinite ones), two monodromy
transformations and thus two zeta-functions $\zeta^0_\varphi(t)$ and
$\zeta^\infty_\varphi(t).$ For a value $c$ different from $0$ and from $\infty$
one can define the same objects (in particular the zeta-function
$\zeta_\varphi^c(t)$) as zero ones for the germ $\varphi-c=\frac{F-c\,G}{G}$.

\remark{Remark 5} One can easily see that these notions are invariant with
respect to projective transformations of the projective line $\hi.$ Thus
$\zeta_\varphi^c(t)$ has to be defined only for one $c$, say, for $c=0.$
\endremark
\vskip 6pt

For the aim of convenience, in \cite{5} we considered only germs of meromorphic
functions $\varphi=\frac{F}{G}$ with $F(0)=G(0)=0.$ At a point which does not
belong to the indeterminacy set of the germ
$\varphi=\frac{F}{G}$ (i.e., if
$F(0)$ or $G(0)$ is different from $0$) one can use the following version of
that
definition. The Milnor fibre of the germ of the \mf $\varphi$ corresponding to a
value $c\in\hi$ is empty unless $\varphi(0)=c$ (and thus the corresponding
zeta-function $\zeta_\varphi^c(t)$ is equal to $1$). The Milnor fibre of the
germ $\varphi$ corresponding to the value $\varphi(0)$ coincides with the usual
Milnor fibre of the holomorphic germ obtained from $\varphi$ by a projective
transformation of $\hi$ which sends $\varphi(0)$ to $0\in\bc^1\subset\hi$ (and
thus the zeta-function $\zeta_\varphi^{\varphi(0)}(t)$ coincides with the usual
zeta-function of this holomorphic germ).

Let $c$ be a point of the projective line $\hi$. For a point $x\in M,$ let
$\zeta_{f,x}^c(t)$ be the corresponding zeta-function of the germ of
the \mf $f$ at the point $x,$ let $\chi_{f,x}^c$ be its degree $\deg
\zeta_{f,x}^c(t).$

\proclaim{Theorem 2}
$$\zeta_f^c(t)=\int_{\{P=Q=0\}\cup F_c}
\zeta_{f,x}^c(t)\,d\chi,\eqno(1)$$
$$\chi(F_{gen})-\chi(F_c)=\int_{F_c}(\chi_{f,x}^c-1)\,d\chi+
\int_{\{P=Q=0\}}\chi_{f,x}^c\,d\chi.\eqno(2)$$
\endproclaim

\demo{Proof} The proof follows the lines of the proof of Theorem 1 in
\cite{4}. Without any loss of generality one can suppose that $c=0.$
There exists a modification $\pi:X\to M$ of the manifold $M$ which is
an isomorphism outside the set
$\{P=Q=0\}\cup\{f=0\}\cup\{f=\infty\}=\{P=0\}\cup\{Q=0\}$ such that
$\dd=\pi^{-1}(\{P=0\}\cup\{Q=0\})$ is a normal crossing divisor in the
manifold $X.$ Then at each point of the exceptional divisor $\dd$ in a local
system of coordinates one has $P\circ \pi=u\cdot
y_1^{k_1}\cdot\ldots\cdot y_n^{k_n}$,
$Q\circ \pi=v\cdot y_1^{\ell_1}\cdot\ldots\cdot y_n^{\ell_n}$ with $u(0)\not=0$,
$v(0)\not=0$, $k_i\geq 0$ and $\ell_i\geq 0.$ There exist  Whitney
stratifications $S'$ and $S^*$ of $M$ and $X$ respectively such that:
\roster
\item the map $\pi$ is a stratified morphism with respect to these
stratifications;
\item for each stratum $\Xi^*\in S^*$ the germs of the liftings $P\circ
\pi$ and $Q\circ\pi$ at points of $\Xi^*$ have normal forms
$u\cdot y_1^{k_1}\cdot\ldots\cdot y_n^{k_n}$ and
$v\cdot y_1^{\ell_1}\cdot\ldots\cdot y_n^{\ell_n}$
 where
$(k_1,\ldots,k_n)$ and $(\ell_1,\ldots,\ell_n)$
do not depend on a point of $\Xi^*.$
\endroster

One applies the following version of the
formula of A'Campo (\cite{1}) and also its local variant.
Let $S_{k,\ell}$ be the set of
points of
$X$ in a neighbourhood of which  the functions $P\circ \pi$ and
$Q\circ \pi$
in some local coordinates have the forms $u\cdot y_1^k$ and $v\cdot y_1^\ell$
respectively
($u(0)\not=0$, $v(0)\not=0$).

\proclaim{Statement 2}
$$\zeta_f^0(t)=\prod_{k>\ell\geq 0}(1-t^{k-\ell})^{\chi(S_{k,\ell})}.$$
\endproclaim

After that, the arguments from Theorem 1 in \cite{4} work literally.
The difference between
$(\chi_{f,x}^c-1)$ and $\chi_{f,x}^c$
in the two integrals in $(2)$ reflects the fact that the Euler characteristic
of the local level set $F_c\cap B_\varepsilon (x)$ ($B_\varepsilon (x)$
is the ball of small radius $\varepsilon$ centred at the point $x$) of the
germ
$f$ is equal to $1$ at a point $x$ of the level set $F_c$ and is equal to
$0$ at a point
$x$ of the indeterminacy set $\{P=Q=0\}$. In the first case
this local level set is contractible and in the second one it is the
difference between two contractible sets. $\square$
\enddemo

Let us denote $(-1)^{n-1}$ times the first integral in $(2)$ by $\mu_f(c)$ and
$(-1)^{n-1}$ times the second one by $\lambda_f(c).$
Let
$\mu_f=\sum\limits_{c\in\hi} \mu_f(c)$ and
$\lambda_f=\sum\limits_{c\in\hi}
\lambda_f(c)$
(in each sum only finite number of summands are different from zero).

\proclaim{Theorem 3}
$$\mu_f+\lambda_f=(-1)^{n-1}\left(2\cdot\chi(F_{gen})
-\chi(M)+\chi(\{P=Q=0\})\right).$$
\endproclaim

\demo{Proof} One has
$$\int_{\hi}
\chi(F_c)\,d\chi=\chi(M\setminus\{P=Q=0\})=\chi(M)-\chi(\{P=Q=0\}).$$
Therefore
$$\chi(M)-\chi(\{P=Q=0\})=\int_{\hi} \chi(F_{gen})\,d\chi+\int_{\hi}
(\chi(F_c)-\chi(F_{gen}))\,d\chi=$$
$$=2\cdot\chi(F_{gen})-(-1)^{n-1}
\sum_{c\in\hi}(\mu_f(c)+\lambda_f(c))=2\cdot\chi(F_{gen})+(-1)^n
(\mu_f+\lambda_f).\qquad \square$$
\enddemo

Let $\widetilde f$ be the restriction of $f$ to $M\setminus\{Q=0\}$,
$\widetilde f: M\setminus\{Q=0\}\to \bc=\hi\setminus\{\infty\}$.
Notice that the fibres of both maps $f$ and $\widetilde f$ 
over values $c\in\bc$ coincide.

\proclaim{Corollary 1}
$$\chi(F_{gen})=\chi(M)-\chi(\{Q=0\})+(-1)^{n-1}(\lambda_f-
\lambda_f(\infty)+\mu_f-\mu_f(\infty)).$$
\endproclaim

Let $f$ be the \mf on the
projective space $\hip$ defined by a polynomial $P$ in $n$ variables (see
below). If $P$ has only isolated critical points in $\bc^n$, then $\mu_f(c)$ is
the sum of the Milnor numbers of the critical points of the polynomial $P$ with
critical value $c$, $\lambda_f(c)$ is equal to the invariant $\lambda_P(c)$
studied in \cite{3}. Therefore $\mu_f(c)$ and $\lambda_f(c)$ can be
considered as generalizations of those invariants. One has
$\mu_f=\mu_P+\mu_f(\infty)$,
$\lambda_f=\lambda_P+\lambda_f(\infty),$ where
$\mu_P=\sum\limits_{c\in\bc} \mu_P(c)$,
$\lambda_P=\sum\limits_{c\in\bc}
\lambda_P(c).$ 
Notice that in this case Corollary 1 turns into 
the well known formula
$\chi(F_{gen})=1+(-1)^{n-1}(\lambda_P+\mu_P).$

\bigbreak
\finparrafo

\head \sec Applications to polynomials
\endhead

 A polynomial $P:\bc^n\to\bc$ defines a \mf
$f=\frac{\widetilde P}{x_0^d}$ on the projective space
$\hip$ ($d=\deg P$). For any $c\in\hi$ the local monodromy of the
polynomial and its zeta-function $\zeta_P^c(t)$ are defined (in fact they
coincide with those of the \mf $f$). The
described technique gives the following statements for polynomials.
Let us remind that for $x\in\{P=c\}\subset\bc^n$ the zeta-function
$\zeta_{f,x}^c(t)$ is the usual zeta-function $\zeta_{P,x}^c(t)$ of the germ
of the polynomial $P$ at $x.$

\proclaim{Theorem 4} For $c\in\bc\subset\hi,$
$$\zeta_P^c(t)=\left( \int_{\{{\widetilde P}=0\}\cap \bc\bp^{n-1}_\infty}
\zeta_{f,x}^c(t)\,d\chi\right)\cdot\left(\int_{\{{P}=c\}}
\zeta_{P,x}^c(t)\,d\chi\right).\eqno(3)$$
For the infinity value,
$$\zeta_P^\infty(t)=\int_{\bc\bp^{n-1}_\infty}
\zeta_{f,x}^\infty(t)\,d\chi.$$
For a generic $c'\in \bc$,
$$\chi(\{P=c'\})-\chi(\{P=c\})=\int_{\{{\widetilde P}=0\}\cap
\bc\bp^{n-1}_\infty}\chi_{f,x}^c\,d\chi+
\int_{\{P=c\}}(\chi_{P,x}^c-1)\,d\chi.$$
\endproclaim

\remark{Remark 6} In \cite{6} we consider the zeta-function
of the local monodromy (corresponding to a finite value $c$) of the polynomial
$P$ near infinity which is just the first factor in the formula $(3)$. If
that zeta-function is different from $1$ then the value $c$ is atypical at
infinity.
\endremark
\vskip 6pt

Let $H=\{R=0\}$ be a hypersurface in $\bc^n$ ($R:\bc^n\to\bc$ is a
polynomial). The polynomial $P$ restricted to the complement of the
hypersurface $H$ defines a $C^\infty$ locally trivial fibration
outside a finite set in $\bc$. For each $c\in\bc$
as well as for $c=\infty$, the local monodromy of this fibration and its
zeta-function
$\zeta_{P/H}^c(t)$ are defined. The described fibration is nothing else but
the fibration for the \mf ${\widetilde f}=\frac{P\cdot R}{R}.$ It implies the
following result.

\proclaim{Theorem 5} For $c\in\bc\subset\hi$,
$$\zeta_{P/H}^c(t)=\int_{\{P=c\}\cup(\{{\widetilde
P}=0\}\cap\bc\bp_\infty^{n-1})}\zeta_{{\widetilde f},x}^c(t)\,d\chi.$$
For  a generic $c'\in\bc$,
$$\chi(F_{c'}\setminus H)-\chi(F_{c}\setminus H)=\int_{\{{\widetilde
P}=0\}\cap
\bc\bp^{n-1}_\infty}\chi_{{\widetilde f},x}^c\,d\chi+\int_{\{P=c\}\cap H}
\chi_{{\widetilde f},x}^c\,d\chi+
\int_{F_{c}\setminus H}(\chi_{P,x}^c-1)\,d\chi.$$
\endproclaim

\bigbreak
\finparrafo


\enddocument